\documentclass[12pt]{article}

\usepackage{amssymb,amsthm,amscd,array}

\let\ssection=\section
\renewcommand{\section}{\setcounter{equation}{0}\ssection}

\newcommand{\bbR}{\mathbb{R}}
\newcommand{\bbRP}{\mathbb{RP}}

\newcommand{\bbC}{\mathbb{C}}

\newcommand{\bbZ}{\mathbb{Z}}

\newcommand{\diag}{\mathrm{diag}}
\newcommand{\cD}{{\mathcal{D}}}

\newcommand{\cE}{{\mathcal{E}}}

\newcommand{\Diff}{\mathrm{Diff}}

\newcommand{\rD}{\mathrm{D}}

\newcommand{\rE}{\mathrm{E}}

\newcommand{\cF}{{\mathcal{F}}}

\newcommand{\rg}{\mathrm{g}}

\newcommand{\cI}{{\mathcal{I}}}
\newcommand{\Id}{\mathrm{Id}}

\newcommand{\Pol}{\mathrm{Pol}}
\newcommand{\cQ}{{\mathcal{Q}}}

\newcommand{\cS}{{\mathcal{S}}}

\newcommand{\SL}{\mathrm{SL}}
\newcommand{\Sl}{\mathrm{sl}}
\newcommand{\SO}{\mathrm{SO}}

\newcommand{\so}{\mathrm{o}}

\newcommand{\ssp}{\mathrm{sp}}

\newcommand{\Symm}{\mathrm{Symm}}

\newcommand{\Vect}{\mathrm{Vect}}

\newcommand{\Vol}{\mathrm{Vol}}

\newcommand{\half}{\frac{1}{2}}
\newcommand{\fg}{\mathfrak{g}}

\def\mod#1{\left|{#1}\right|}

\begin{document}



\def\a{\alpha}
\def\b{\beta}
\def\d{\delta}
\def\g{\gamma}
\def\om{\omega}
\def\r{\rho}
\def\s{\sigma}
\def\vfi{\varphi}
\def\vr{\varrho}
\def\l{\lambda}
\def\m{\mu}
\def\implies{\Rightarrow}

\oddsidemargin .1truein
\newtheorem{thm}{Theorem}[section]
\newtheorem{lem}[thm]{Lemma}
\newtheorem{cor}[thm]{Corollary}
\newtheorem{pro}[thm]{Proposition}
\newtheorem{ex}[thm]{Example}
\newtheorem{rmk}[thm]{Remark}
\newtheorem{defi}[thm]{Definition}


\title{Methods of Equivariant Quantization}

\author{C.~Duval\thanks{
Universit\'e de la M\'editerran\'ee and CPT-CNRS,
Luminy Case 907, F--13288 Marseille, Cedex~9, FRANCE;
mailto:duval@cpt.univ-mrs.fr
}
\and
P.~Lecomte\thanks{
Institut de Math\'ematiques, Universit\'e de Li\`ege, Sart Tilman, Gde
Traverse, 12 (B~37), B--4000 Li\`ege, BELGIUM;
mailto:plecomte@ulg.ac.be
}
\and
V.~Ovsienko\thanks{
CNRS, Centre de Physique Th\'eorique,
CPT-CNRS, Luminy Case 907,
F--13288 Marseille, Cedex~9, FRANCE;
mailto:ovsienko@cpt.univ-mrs.fr
}
}

\date{}

\maketitle

\bigskip

\thispagestyle{empty}

\begin{abstract}
This article is a survey of recent work \cite{LO,DO,DLO,Lec1} 
developing a new approach to quantization based on the
equivariance with respect to some Lie group of symmetries.
Examples are provided by conformal and projective
differential geometry: given a smooth manifold $M$ endowed with
a flat conformal/projective structure, we establish a 
canonical isomorphism between the space of symmetric contravariant
tensor fields on~$M$ 
and the space of differential operators on~$M$.
This leads to a notion of conformally/projectively
invariant star-product on $T^*M$.
\end{abstract}

\vskip1cm
\noindent
\textbf{Keywords:} Quantization, Lie algebras of vector fields, 
hypergeometric functions
 star-products.

\vskip1cm
\noindent
Proc. of the workshop ``Noncommutative Differential Geometry and its Applications to Physics'',
Shonan-Kokusaimura, Japan, May 31-June 4, 1999.

\newpage

\section{Introduction}

This paper is a survey of recent articles
\cite{LO,DO,DLO,Lec1} (see also \cite{CMZ})
developing a new viewpoint on quantization.

\medskip
\noindent
{\bf 1.1}
Let $\cD(M)$ be the space of differential operators on a smooth manifold $M$ and
$\cS(M)(=\Pol(T^*M))$
the space of symbols, that is, of smooth functions on $T^*M$
polynomial on fibers. Note that $\cS(M)$ is naturally isomorphic
to the space $\Gamma(STM)$ of symmetric contravariant tensor fields on $M$.

A \textit{quantization map} (inverse to a \textit{symbol map}) is
a linear map 
\begin{equation}
\label{QuMap}
\cQ:\cS(M)\to\cD(M)
\end{equation}
satisfying the following natural
condition: for every $P\in\cS(M)$
the principal symbol of $\cQ(P)$ coincides with the homogeneous higher-order
part of $P$.

We will be considering a Lie group $G$ acting on $M$ and impose the
$G$-equivariance condition; we will then seek a map (\ref{QuMap})
commuting with the $G$-action.
We will be, furthermore, interested in the situation where the equivariance
condition determines the quantization map \textit{in a unique fashion}.

\medskip
\noindent
{\bf 1.2}
Our main examples of symmetry groups will be:
\begin{equation}
\label{SyGr}
G=\SL(n+1,\bbR)
\quad
\hbox{and}
\quad
G=\SO(p+1,q+1),
\quad
p+q=n
\quad
(=\dim{M})
\end{equation}
acting in the standard way on $S^n\,(\to\bbRP^n)$
and on $S^p\times{}S^q\;(\to (S^P\times S^q)/\bbZ_2)$ (respectively).
These two cases are related to 
\textit{projective} and \textit{conformal} differential geometry.
Our main purpose is to introduce the projectively/conformally equivariant
quantization and prove its uniqueness.

\medskip
\noindent
{\bf 1.3}
The crucial property common to these two cases is \textit{maximality}
in the following sense. The corresponding Lie algebras:
\begin{equation}
\label{SyAl}
\fg=\Sl(n+1,\bbR)
\quad
\hbox{and}
\quad
\fg=\so(p+1,q+1),
\quad
p+q\geq3
\end{equation}
can be both viewed as maximal subalgebras of the Lie algebra 
$\Vect_{\Pol}(\bbR^n)$ of all polynomial vector fields on $\bbR^n$
(cf. \cite{LO,BL}). We believe
(although this remains unproven) that the uniqueness of the equivariant
quantization map is due to
the maximality of the algebras of symmetry. 

\goodbreak

\medskip
\noindent
{\bf 1.4}
Consider the group of all diffeomorphisms $\Diff(M)$ and the Lie algebra
of all smooth vector fields $\Vect(M)$.
An important aspect of our approach is the study of $\Diff(M)$-~and
$\Vect(M)$-module structure on the space of differential operators .

Let~$\cF_\l$ be the space of \textit{tensor densities} of degree $\l$ on $M$:
\begin{equation}
\label{TenDen}
\varphi=
f(x^1,\ldots,x^n)
\left|
dx^1\wedge\cdots\wedge{dx^n}
\right|^\l
\end{equation}
that is, of sections of the line bundle
$\Delta_\l(M)=\mod{\Lambda^nT^*M}^{\otimes\l}$ over $M$.
Denote $\cD_{\l,\m}(M)$ the space
of differential operators
\begin{equation}
\label{DofOp}
A:\cF_\l\to\cF_\m.
\end{equation}
This space is naturally a module over $\Diff(M)$ and $\Vect(M)$.
Therefore, for a Lie subgroup $G\subset\Diff(M)$, one obtains a
$G$-action on $\cD_{\l,\m}(M)$.

\begin{rmk}
{\rm
Assume that $M$ is orientable and choose a volume form on $M$, one then
naturally identifies $\cD_{\l,\m}(M)$ with $\cD(M)$ (i.e. with $\cD_{0,0}(M)$)
as a vector space. One can, therefore, consider
$\cD_{\l,\m}(M)$ as a two-parameter family of $\Diff(M)$- and 
$\Vect(M)$-actions on the same space, $\cD(M)$.
}
\end{rmk}

It worth noticing that many important examples of differential operators
acting on tensor densities have already been studied by the classics.
The most important case $\cD_{\half,\half}$ 
(of operators on \textit{half}-densities) naturally arises in geometric
quantization (see e.g. \cite{Kir}).
We refer to recent articles \cite{DO1,LMT,GO,Mat} 
(and references therein) for a systematic study
and classification of the modules $\cD_{\l,\m}(M)$.

\medskip
\noindent
{\bf 1.5}
In many cases we have to assume that the $G$-action cannot be
defined \textit{globally} on~$M$. An example is provided 
when $M$ is not compact (e.g., $M=\bbR^n$) and one can only deal
with a $\fg$-action. Moreover, there are many cases when even
the $\fg$-action is only locally-defined.
The most general framework is provided by the notion of a flat $G$-\textit{structure}
(defined by a local action of $G$ on
$M$, compatible with a local identification of~$M$ with some homogeneous space~$G/H$).
Existence of such structures on a given smooth manifold is a difficult
(and widely open) problem of topology (see e.g. \cite{Thu}). 

We will assume that $M$ is endowed either 
with a flat projective or a flat conformal structure.

\medskip
\noindent
{\bf 1.6}
In the particular case $n=1$, both conformal and projective structures coincide. We
refer to \cite{CMZ} (ee also \cite{Wil}) for a thorough study 
of $\Sl(2,\bbR)$-equivariant
quantization and of the corresponding invariant star-product.

\section{Existence results}\label{Res}

In this section we formulate our most general existence theorems.
We will treat the basic examples and give some explicit formul{\ae}
later.

\subsection{Symbols with values in tensor densities}\label{SymbDens}

The space $\cS(M)$ of contravariant symmetric tensor
fields on~$M$ is naturally a $\Diff(M)$- and $\Vect(M)$-module. 
We will, however, need to consider a one-parameter family of 
$\Diff(M)$- and $\Vect(M)$-actions on this space.
Put
\begin{equation}
\cS_\d
=
\cS(M)
\otimes_{C^\infty(M)}
\cF_\d.
\label{P}
\end{equation}
The space $\cS_\d$ is also, naturally, a $\Diff(M)$- and $\Vect(M)$-module.

\subsection{Action of $\Vect(M)$ on the space of symbols}

Space $\cS(M)$ is a Poisson algebra isomorphic as a $\Vect(M)$-module to
the space of smooth functions on $T^*M$ polynomial on the fibers.
The action of $X\in\Vect(M)$ on $\cS$ is given by the Hamiltonian vector field
\begin{equation}
\label{Hamilton}
L_X=
\frac{\partial{}X}{\partial\xi_i}\,
\frac{\partial}{\partial{}x^i}
-
\frac{\partial{}X}{\partial{}x^i}\,
\frac{\partial}{\partial\xi_i}\,,
\end{equation}
where $(x^i,\xi_i)$ are local coordinates on $T^*M$ 
(we identified vector fields on~$M$
with the first-order polynomials on $T^*M$, that is
$X=X^i\xi_i$ in (\ref{Hamilton}))\footnote{
Throughout this paper the summation over repeated indices is understood.}.
The Poisson bracket on $\cS(M)$ is usually called the
(symmetric) Schouten bracket (see e.g.
\cite{Fuc}).

The $\Vect(M)$-action on $\cS_\d$ is of the form:
\begin{equation}
\label{SymbAct}
L^\d_X=
L_X
+\d{}\rD(X)
\end{equation}
where
\begin{equation}
\label{DivEq}
\rD(X)=
\partial_iX^i,
\qquad
\partial_i=
\partial/\partial{x^i}.
\end{equation}
Note that the formula (\ref{SymbAct}) does not depend on the choice
of local coordinates.
The space $\cS_\d$ is naturally the space of symbols corresponding to the
space of differential operators $\cD_{\l,\m}$ such that $\m-\l=\d$.

\subsection{Isomorphism of $\fg$-modules}

Let now $M$ be a endowed with a flat projective or conformal structure
(see e.g. \cite{Thu} and Section \ref{Coord} for precise definition).
The most general results of \cite{LO,Lec1,DO,DLO} can be formulated as follows.

\begin{thm}
\label{IsomGen1}
For generic values of $\d$ there exists an isomorphism of
$\Sl(n+1,\bbR)$- or $\so(p+1,q+1)$-modules (respectively):
\begin{equation}
\widetilde{\cQ}_{\l,\m}:\cS_\d
\stackrel{\cong}\rightarrow
\cD_{\l,\m},
\qquad\d=\m-\l
\label{prequant}
\end{equation}
which is unique provided the principal symbol be preserved at each order.
\end{thm}

There exist values of $\d$ such that the
the modules $\cS_\d$ and $\cD_{\l,\m}$ with $\m-\l=\d$
are not isomorphic for generic $\l$.
These particular values are called \textit{resonant}.

\begin{rmk}
{\rm
\hfill\break\indent
(a)
In the projective case the
resonance values of $\d$ are $\d=k/(n+1)$, where $k$ is integer $\geq{}n+1$
(see \cite{Lec1}).

(b)
In the conformal case the structure of resonances is much more complicated
(see \cite{DLO}).

(c)
In the one-dimensional case, $M=S^1$ or $\bbR$, the above theorem still holds true but the
resonances are simply $\d=1,\frac{3}{2},2,\frac{5}{2},\ldots$ and appear in
\cite{CMZ,Gar}. (The projective or conformal structure is then replaced
by the natural $\Sl(2,\bbR)$ action.)
}
\end{rmk}

\subsection{Definition of quantization map}

Let us introduce a new operator on symbols that will eventually insure the
symmetry of the corresponding differential operators. 
Define
$\cI_\hbar:\cS_\d\to\cS_\d[[i\hbar]]$ by
\begin{equation}
\cI_\hbar(P)(\xi)=P(i\hbar\,\xi).
\label{Ihbar}
\end{equation}
Note that we will understand $\hbar$ either as a formal parameter or as a fixed real
number, depending upon the context.
We will introduce the quantization map
$\cQ_{\l,\m;\hbar}:\cS_\d\to\cD_\l[[i\hbar]]$ depending on $\hbar$:
\begin{equation}
\cQ_{\l,\m;\hbar}=\widetilde{\cQ}_{\l,\m}\circ\cI_\hbar
\label{quantGen}
\end{equation}
where $\widetilde{\cQ}_{\l,\m}$ is the isomorphism from Theorem \ref{IsomGen1}.

\subsection{Symmetric quantized operators}

Let us recall that if $\l+\m=1$, there exists, for compactly-supported
densities, a $\Vect(M)$-invariant pairing $\cF^{\bbC}_\l\otimes\cF^{\bbC}_\m\to\bbC$ defined by
\begin{equation}
\varphi\otimes\psi\mapsto\int_M{\!\overline{\varphi}\,\psi}.
\label{pairing}
\end{equation}

We can now formulate the important
\begin{cor}
For generic $\d$ and $\l+\m=1$ the quantization
\begin{equation}
\check{P}
=
\cQ_{\frac{1-\d}{2},\frac{1+\d}{2};\hbar}(P)
\label{quantGenSymm}
\end{equation}
of any symbol $P\in\cS_\d$ is a symmetric (formally self-adjoint)
operator.
\end{cor}
\begin{proof}
Let us denote by $A^*$ the adjoint of $A\in\cD_{\l,1-\l}$ with respect to the
pairing~(\ref{pairing}). Consider the symmetric operator
$$
\Symm(\cQ_{\l,\m;\hbar}(P))
=
\half\Big(
\cQ_{\l,\m;\hbar}(P)+(\cQ_{\l,\m;\hbar}(P))^*
\Big)
$$
which exists whenever $\l+\m=1$. Notice that it has the same principal symbol as
$\cQ_{\l,\m;\hbar}(P)$. Now, the map $\Symm(\cQ_{\l,\m;\hbar})$ is obviously
$\fg$-equivariant. The result follows then from the uniqueness of
$\fg$-module isomorphism ---
see Theorem \ref{IsomGen1}.
\end{proof}

\goodbreak

\subsection{Invariant star-products}

The special and most important value $\d=0$ is non-resonant.
The existence and uniqueness of the isomorphism
(\ref{prequant}) in this case has been proven in \cite{LO,DLO}.

Define an associative bilinear operation (depending on $\l$)
\begin{equation}
*_{\l;\hbar}:\cS\otimes\cS\to\cS[[i\hbar]]
\label{starOperation}
\end{equation}
such that
\begin{equation}
\cQ_{\l;\hbar}(P*_{\l;\hbar}Q)
=
\cQ_{\l;\hbar}(P)\circ\cQ_{\l;\hbar}(Q).
\label{star}
\end{equation}

Recall that an associative operation $*_\hbar:\cS\otimes\cS\to\cS[[i\hbar]]$ is called
a \textit{star-product} if it is of the form
\begin{equation}
P*_\hbar{}Q=PQ+\frac{i\hbar}{2}\{P,Q\}+O(\hbar^2)
\label{starDef}
\end{equation}
where $\{\cdot,\cdot\}$ stands for the Poisson bracket on $T^*M$, and is given by
bi-differential operators at each order in~$\hbar$.

\begin{thm}\label{starProduct}
The associative, conformally invariant, operation $*_{\l;\hbar}$ defined by (\ref{star})
is a star-product if and only if $\l=\half$.
\end{thm}

Let us emphasize that this theorem provides us precisely with the value of
$\l$ used in geometric quantization and, in some sense, links the latter to
deformation quantization.

\subsection{Quantum Hamiltonians in (pseudo-)Riemannian case}
\label{quantumHamiltonian}

Assume now that a (pseudo-)Riemannian metric $\rg$ on $M$ is fixed.

To recover the traditional Schr\"odinger picture of quantum mechanics,
one needs to associate to the operator $\check{P}$ resulting from our quantization
map (\ref{quantGenSymm}) an operator 
\begin{equation}
\hat{P}:\cF^{\bbC}_0\to\cF^{\bbC}_0
\label{quantumOp}
\end{equation}
on the space of complex-valued functions on a conformally flat
manifold $(M,\rg)$.

\begin{defi}
{\rm 
Using the natural identification $\cF^{\bbC}_0\to\cF^{\bbC}_\l$ between tensor densities and smooth
functions given by 
\begin{equation}
f\mapsto{}f\,\mod{\Vol_\rg}^\l,
\label{identification}
\end{equation}
one can define the differential operator $\hat{P}$ by the commutative
diagram
\begin{equation}
\begin{CD}
\cF^{\bbC}_0 @> \textstyle{\hat{P}} >> \cF^{\bbC}_0 \strut\\
@V{\mod{\Vol_\rg}^\l}VV @VV{\mod{\Vol_\rg}^\m}V \strut\\
\cF^{\bbC}_\l @> \textstyle{\check{P}} >> \cF^{\bbC}_\m \strut
\end{CD}
\label{TheNewDiagram}
\end{equation}
where $\check{P}$ is given by (\ref{quantGenSymm}) in the case $\l+\m=1$.
}
\end{defi}
We call the operator $\hat{P}$ corresponding to a polynomial $P$ the
quantum Hamiltonian.

\begin{rmk}
{\rm 
An important feature of equivariant quantization in the 
conformal case is, of course, that the metric $\rg$ can be chosen
conformally flat (as a representative of the conformal structure on $M$).
}
\end{rmk}

\goodbreak

\section{Examples: quantizing quadratic Hamiltonians}\label{SecEx}

Let us consider a homogeneous second order polynomial $H$ on $T^*M$.
In arbitrary local coordinates it is of the form:
\begin{equation}
H=
\rg^{ij}(x)\xi_i\xi_j.
\label{HomHaM}
\end{equation}
To illustrate the general Theorem \ref{IsomGen1} we will give explicit
formul{\ae} for the quantization in this case.

We will need special coordinate systems adopted to projective and conformal structures.

\subsection{Local coordinates}\label{Coord}

The following description can be considered as a definition of
flat projective and conformal structures.

\begin{defi}
{\rm
A smooth manifold $M$ is endowed with a 
projective (or conformal) structure if
there exists an atlas $U_i$ on $M$ with local coordinates
$(x^1,\ldots,x^n)_i$ in each chart $U_i$ that are accorded on 
$U_i{\cap}U_j$ in the following way. 

\goodbreak

(1) \textit{Projective case}:
The linear span of the following $n(n+2)$ vector fields
\begin{equation} 
\label{sl}
\left\langle
\partial_i\,, 
\quad 
x^i\partial_j\,, 
\quad
x^i{\cal E}\,
\right\rangle
\;,
\qquad
\hbox{where}
\quad
{\cal E} = x^j \partial_j\,.
\end{equation}
is independent on the choice of a chart 
(i.e., coincide on each $U_i{\cap}U_j$). 
The vector fields (\ref{sl}) obviously generate the Lie algebra
isomoprhic to $\Sl(n+1,\bbR)$ and, therefore, one obtains a (locally defined)
$\Sl(n+1,\bbR)$-action on $M$.

It worth noticing that coordinate changes on $U_i{\cap}U_j$ in this case
are nothing but linear-fractional transformations.

(2) \textit{Conformal case}:
The linear span of the following $(n+2)(n+1)/2$ vector fields
\begin{equation}
\left\langle
\partial_i,\;
x_i\partial_j-
x_j\partial_i,\;
\cE,\;
x_jx^j\partial_i-
2x_i\cE
\right\rangle
\label{confGenerators}
\end{equation}
coincide on each $U_i{\cap}U_j$
(we have used the notation $x_i=g_{ij}x^j$ where the flat
metric $g=\diag(1,\ldots,1,-1,\ldots,-1)$ is of signature $(p,q)$). 
Note that the vector fields~(\ref{confGenerators}) generate the Lie algebra
isomorphic to $\so(p+1,q+1)$.

The coordinate changes are then given by the conformal tranformations.
}
\end{defi}

\subsection{The quantization map in the second-order case}
\label{SecQuan}

We are now ready to give explicit expressions of our quantization maps.
For the sake of simplicity we restrict ourself to the most interesting
case $\l=\m=1/2$ (see~\cite{LO,DO,DLO} for the general formul{\ae}).
We will use coordinates adopted to a projective or conformal structure
on $M$.

\medskip
\noindent
(1) \textit{Projective case}:
The quantization map (\ref{quantGen})
associates to the polynomial (\ref{HomHaM}) the following second-order 
differential operator
\begin{equation}
\label{QHamPro}
\cQ_{\half;\hbar}(H)=
-\hbar^2
\left(
\rg^{ij}\partial_i\partial_j+
\partial_j(\rg^{ij})\partial_i+
\frac{(n+1)}{4(n+2)}\partial_i\partial_j(\rg^{ij})
\right)
\end{equation}
on the space $\cF^{\bbC}_{1/2}$ of half-densities.

\medskip
\noindent
(2) \textit{Conformal case}:
One has
\begin{equation}
\begin{array}{rcl}
\cQ_{\half;\hbar}(H) &=&
-\hbar^2
\Big(
\rg^{ij}\partial_i\partial_j+
\partial_j(\rg^{ij})\partial_i\\[12pt]
&&+\displaystyle
\frac{n}{4(n+1)}\partial_i\partial_j(\rg^{ij})
+\frac{n}{4(n+1)(n+2)}\partial_i\partial_i(\rg^{jj})
\Big)
\end{array}
\label{Quantization}
\end{equation}

Let us also recall here for comparison
the explicit expression of the most famous and standard
quantization procedure, called the 
\textit{Weyl} or \textit{symmetric} quantization.

\medskip
\noindent
(3) \textit{Weyl quantization}:
In the case of second-order polynomials 
\begin{equation}
\cQ_{{\rm Weyl}}(H) =
-\hbar^2
\left(
\rg^{ij}\partial_i\partial_j+
\partial_j(\rg^{ij})\partial_i+
\frac{1}{4}\partial_i\partial_j(\rg^{ij})
\right)
\label{WQuantiz}
\end{equation}
This formula makes sense globally only in the case $M=\bbR^n$ and 
it is then defined uniquely
(up to normalization) by equivariance with respect to the standard
$\ssp(2n,\bbR)$-action on $T^*M=\bbR^{2n}$.

\begin{rmk}
{\rm
It can be checked that:

(a) the explicit formul{\ae} (\ref{QHamPro},\ref{Quantization})
are independent of the choice
of the adopted coordinate system;

(b) in both cases the differential operator
$\cQ_{\half;\hbar}(H)$ is symmetric.
}
\end{rmk}

\subsection{Quantizing the geodesic flow}\label{SecQuan}

Assume now that there exists (a non-degenerate) metric
\begin{equation}
\label{Hmetric}
\rg=
\rg_{ij}dx^idx^j
\end{equation}
corresponding to the polynomial (\ref{HomHaM}), 
i.e. $H=\rg^{-1}$. 

Let us apply the procedure of Section \ref{quantumHamiltonian} to
obtain the quantum Hamiltonian. It is easy to see that the
fact that the operators (\ref{QHamPro}-\ref{WQuantiz}) are symmetric
implies in each case (1), (2), (3) that the corresponding
quantum Hamiltonian has the general form:
\begin{equation}
\hat{H}
=
-\hbar^2
\left(
\Delta_\rg+U
\right)
\label{GeodFlowGen}
\end{equation}
where $\Delta_\rg$ is the Laplace-Beltrami operator and the function $U$
(the potential) is different in each of the three cases.

The potential can be easily computed in each case (1), (2), (3).
In particular, in the conformal case one can choose 
$\rg$ \textit{conformally flat}.
It has been shown in~\cite{DO} that in such a situation
\begin{equation}
\hat{H}
=
-\hbar^2\left(
\Delta_\rg-\frac{n^2}{4(n-1)(n+2)}\,R_\rg
\right)
\label{QuantGeodFlow}
\end{equation}
where $R_\rg$ stands for the scalar curvature of $(M,\rg)$.
It worth noticing that this result differs from other known
versions of the quantization of the geodesic flow
(cf. \cite{det} and references therein).

\section{Universal formula in the projective case}\label{ConfInvOp}

In the projective case there exists an explicit formula 
for the quantization map that is valid
for polynomial functions on $T^*M$ of an arbitrary order $k$ (see \cite{LO,Lec1}).
In this section we will give a new (but still equivalent) expression
of the projectively equivariant quantization map.
As above, we will consider only the case $\l=\m=1/2$.

\subsection{The Euler and divergence operators}

Fix an arbitrary system of local coordinates $(x^1,\ldots,x^n)$ on $M$ and
introduce the following differential operators acting 
\textit{locally} on the space 
of symbols (i.e., on the space of polynomials 
$\bbC[x^1,\ldots,x^n,\xi_1,\ldots,\xi_n]$):
\begin{equation}
\label{ED}
\rE=
\xi_i\frac{\partial}{\partial \xi_i}+\frac{n}{2},
\qquad
\rD=\frac{\partial}{\partial \xi_i}\frac{\partial}{\partial x^i}
\end{equation}

\goodbreak

\subsection{The confluent hypergeometric function}

Recall that the confluent hypergeometric function is a series in $z$
depending on two parameters $a,b\in{\bbR}$ (or $\bbC$):
$$
F\left(\left.
\begin{array}{c}
a\\b
\end{array}\right|z\right)
=
\sum_{m=0}^\infty
{
\frac{(a)_m}{(b)_m}\,\frac{z^m}{m!}
}
$$
where
$$
(a)_m:=a(a+1)\cdots(a+m-1)
$$
(see, e.g., \cite{GKD}).

\subsection{The explicit formula}

Identify locally polynomials with differential operators 
by the normal ordering map
\begin{equation}
\s:A_k^{{i_1}\ldots{i_k}}
\partial_{i_1}
\cdots
\partial_{i_k}
\mapsto
A_k^{{i_1}\ldots{i_k}}
\xi_{i_1}
\cdots
\xi_{i_k}
\label{sigma}
\end{equation}
Using this identification we will work with differential operators 
as with polynomials.

The following theorem provides an explicit expression
for the projectively equivariant quantization map (\ref{quantGen})
in terms of a ``non-commutative'' confluent function.

\begin{thm}
\label{ExpThm}
The projectively equivariant quantization is of the form
\begin{equation}
\cQ_{\half;\hbar}=F\left(\left.
\begin{array}{r}
2\rE\\\rE
\end{array}\right|\frac{i\hbar\rD}{4}\right)
\label{QPE}
\end{equation}
for every system of local coordinates adopted to the projective structure.
\end{thm}

In other words, one has
\begin{equation}
\cQ_{\half;\hbar}
=
\sum_{m=0}^\infty{
C_m(\rE)\,(i\hbar\rD)^m
}
\label{Ansatz}
\end{equation}
where the expression
\begin{equation}
C_m(\rE)
=
\frac{1}{m!}\,
\frac{(\rE+1/2)(\rE+3/2)\cdots(\rE+m-1/2)}
{(2\rE+m)(2\rE+m+1)\cdots(2\rE+2m-1)}
\label{Coeff}
\end{equation}
is also understood as an infinite (Taylor) series in $\rE$.

\goodbreak

\begin{rmk}
\hfill\break\indent
{\rm
(a)
If the manifold $M$ is endowed with a projective structure,
then choosing an adapted coordinate system,
the formula (\ref{QPE}) is well-defined globally on $T^*M$.

(b) 
One easily checks that for the second-order polynomials
the formula (\ref{QPE}) coincides with (\ref{QHamPro}).
}
\end{rmk}

The expression (\ref{QPE}) is a particular case (for $\l=1/2$)
of the quantization 
map given in \cite{LO} after restriction to the space of $k$-th order
polynomials.
Unfortunately, we do not know the analog of the
formula (\ref{QPE}) in the conformal case.

\subsection{Comparison with the Weyl quantization}\label{ComparWey}

The Weyl quantization map, $\cQ_{\rm Weyl}$, retains the very elegant form:
\begin{equation}
\begin{array}{rcl}
\displaystyle
\cQ_{\rm Weyl} 
&=&
\displaystyle\exp\Big(\frac{i\hbar}{\!2}\rD\Big)\\[10pt] &=&
\displaystyle
\Id
+\frac{i\hbar}{\!2}\rD
-\frac{\hbar^2}{\!8}\rD^2
+
O(\hbar^3)
\end{array}
\label{WeylMap}
\end{equation}
(see \cite{AW}).
Comparing (\ref{WeylMap}) with (\ref{QPE}) one notes that the
exponential is, indeed, the most elementary hypergeometric function.

\subsection{An open problem}

It is well-known that the Weyl quantization can be extended for a much
larger class of functions on $\bbR^{2n}$ than polynomials, namely to
the space of pseudo-differential symbols.

We formulate here the following problem.
It it true that the formula (\ref{QPE}) is also valid for pseudo-differential
symbols~? If the answer is positive,
it would be interesting to compute its integral expression.

\newpage


\end{document}